\documentclass[a4paper,12pt]{amsart}
\usepackage{amsfonts}
\usepackage{amssymb}
\usepackage{ifthen}
\usepackage{graphicx}
\nonstopmode \numberwithin{equation}{section}
\usepackage[usenames]{color}
\newtheorem{thm}{Theorem}
\newtheorem{lem}{Lemma}
\newtheorem{cor}{Corollary}
\newtheorem{prop}{Proposition}
\newtheorem{ques}{Question}

\newtheorem{cl}{Claim}
\newtheorem{ca}{Case}
\newtheorem{sca}{Subcase}
\newtheorem{scl}{Subclaim}
\newtheorem{conj}[equation]{Conjecture}

\theoremstyle{definition}
\newtheorem{defn}{Definition}

\newtheorem{op}[equation]{Open Problem}
\newtheorem{rem}{Remark}
\newtheorem{examp}{Example}

\newcounter {own}
\def\theown {\thesection       .\arabic{own}}

\newenvironment{pf}[1][]{%
 \vskip 3mm
 \noindent
 \ifthenelse{\equal{#1}{}}%
  {{\slshape Proof. }}%
  {{\slshape #1.} }%
 }%
{\qed\bigskip}

\newcounter{alphabet}
\newcounter{tmp}
\newenvironment{Thm}[1][]{\refstepcounter{alphabet}%
\bigskip%
\noindent%
{\bf Theorem \Alph{alphabet}}%
\ifthenelse{\equal{#1}{}}{}{ (#1)}%
{\bf .} \itshape}{\vskip 8pt}

\makeatletter
\newcommand{\Ref}[1]{\@ifundefined{r@#1}{}{\setcounter{tmp}{\ref{#1}}\Alph{tmp}}}
\makeatother

\def\be{\begin{equation}}
\def\ee{\end{equation}}

\newcommand{\bee}{\begin{enumerate}}
\newcommand{\eee}{\end{enumerate}}

\newcommand{\blem}{\begin{lem}}
\newcommand{\elem}{\end{lem}}
\newcommand{\bthm}{\begin{thm}}
\newcommand{\ethm}{\end{thm}}
\newcommand{\bcor}{\begin{cor}}
\newcommand{\ecor}{\end{cor}}
\newcommand{\beg}{\begin{exam}}
\newcommand{\eeg}{\end{exam}}
\newcommand{\begs}{\begin{examples}}
\newcommand{\eegs}{\end{examples}}
\newcommand{\bdefe}{\begin{defn}}
\newcommand{\edefe}{\end{defn}}
\newcommand{\bprob}{\begin{prob}}
\newcommand{\eprob}{\end{prob}}
\newcommand{\bques}{\begin{ques}}
\newcommand{\eques}{\end{ques}}
\newcommand{\bei}{\begin{itemize}}
\newcommand{\eei}{\end{itemize}}
\newcommand{\bcon}{\begin{conj}}
\newcommand{\econ}{\end{conj}}
\newcommand{\bop}{\begin{op}}
\newcommand{\eop}{\end{op}}

\newcommand{\bca}{\begin{ca}}
\newcommand{\eca}{\end{ca}}
\newcommand{\bsca}{\begin{sca}}
\newcommand{\esca}{\end{sca}}

\newcommand{\bcl}{\begin{cl}}
\newcommand{\ecl}{\end{cl}}

\newcommand{\bscl}{\begin{scl}}
\newcommand{\escl}{\end{scl}}

\newcommand{\bcons}{\begin{conjs}}
\newcommand{\econs}{\end{conjs}}
\newcommand{\bprop}{\begin{propo}}
\newcommand{\eprop}{\end{propo}}
\newcommand{\br}{\begin{rem}}
\newcommand{\er}{\end{rem}}
\newcommand{\brs}{\begin{rems}}
\newcommand{\ers}{\end{rems}}
\newcommand{\bo}{\begin{obser}}
\newcommand{\eo}{\end{obser}}
\newcommand{\bos}{\begin{obsers}}
\newcommand{\eos}{\end{obsers}}
\newcommand{\bpf}{\begin{pf}}
\newcommand{\epf}{\end{pf}}
\newcommand{\ba}{\begin{array}}
\newcommand{\ea}{\end{array}}
\newcommand{\beq}{\begin{eqnarray}}
\newcommand{\beqq}{\begin{eqnarray*}}
\newcommand{\eeq}{\end{eqnarray}}
\newcommand{\eeqq}{\end{eqnarray*}}

\newcounter{minutes}
\divide\time by 60
\newcounter{hours}
\multiply\time by 60 \addtocounter{minutes}{-\time}
\begin{document}
\bibliographystyle{amsplain}
\title [] {Coefficient estimates and the Fekete-Szeg\H{o}
problem for certain classes of polyharmonic mappings}

\def\thefootnote{}
\footnotetext{ \texttt{\tiny File:~\jobname .tex,
          printed: \number\day-\number\month-\number\year,
          \thehours.\ifnum\theminutes<10{0}\fi\theminutes}
} \makeatletter\def\thefootnote{\@arabic\c@footnote}\makeatother

\author{J. Chen}
\address{J. Chen, Department of Mathematics,
Hunan Normal University, Changsha, Hunan 410081, People's Republic
of China.} \email{jiaolongchen@sina.com}

\author{A. Rasila}
\address{A. Rasila, Department of Mathematics,
Hunan Normal University, Changsha, Hunan 410081, People's Republic
of China, and
Department of Mathematics and Systems Analysis, Aalto University, P. O. Box 11100, FI-00076 Aalto,
 Finland.} \email{antti.rasila@iki.fi}

\author{X. Wang${}^{~\mathbf{*}}$}
\address{X. Wang, Department of Mathematics,
Hunan Normal University, Changsha, Hunan 410081, People's Republic
of China.} \email{xtwang@hunnu.edu.cn}

\keywords{harmonic mapping, polyharmonic mapping, coefficient estimates, the Fekete-Szeg\H{o}
problem.\\
${}^{\mathbf{*}}$ Corresponding author}

\begin{abstract}
We give coefficient estimates for a class of close-to-convex harmonic mappings,
and discuss the Fekete-Szeg\H{o}
problem of it.
We also introduce two classes of polyharmonic
mappings $\mathcal{HS}_{p}$ and $\mathcal{HC}_{p}$, consider the starlikeness and convexity of them, and obtain coefficient estimates on them.
Finally, we give a necessary condition for a mapping $F$ to be in the class $\mathcal{HC}_{p}$.

\end{abstract}

\thanks{The research was partly supported by
NSF of China (No. 11071063) and  Hunan Provincial Innovation
Foundation for Postgraduates (No. 125000-4242). }

\maketitle \pagestyle{myheadings} \markboth{J. Chen, A. Rasila
and X. Wang}{Coefficient estimates and the Fekete-Szeg\H{o}
problem for certain classes of polyharmonic mappings}

\section{Introduction}\label{csw-sec1}

Let $\mathbb{D}$ denote the unit disk $\{z:|z|<1,z\in \mathbb{C}\}$,
 and let $\mathcal{A}$ be the class of functions of the form
\be\label{eq1.1}f(z)=z+\sum_{j=2}^{\infty}a_{j}z^{j} ,\ee
which are analytic in $\mathbb{D}$. Denote by $\mathcal{S}$ the subclass of $\mathcal{A}$
consisting of functions $f\in \mathcal{A}$,
which are univalent.
A continuous mapping $f=u+iv$ is a {\it complex-valued harmonic} mapping in a domain $D\subset \mathbb{C}$
if both $u$ and $v$ are real harmonic in $D$, i.e., $\Delta u=\Delta v=0$, where $\Delta$
is the complex Laplacian operator
$$\Delta=4\frac{\partial^{2}}{\partial
z\partial\bar{z}}:=\frac{\partial^2}{\partial
x^2}+\frac{\partial^2}{\partial y^2}.$$
For a simply connected
domain $D$, we can write $f$ in the
form $f = h+ \overline{g}$, where $h$ and $g$ are analytic (see \cite{du}).
A necessary and sufficient condition for $f$ to be locally
univalent and sense-preserving in $D$ is that $|h'(z)|>|g'(z)|$ for all $z\in D$.

A continuous complex-valued mapping $F$ in $D$ is
{\it biharmonic} if the Laplacian of $F$ is harmonic, i.e.,
$F$ satisfies the equation $\Delta(\Delta F) = 0$. It can be shown
that in a simply connected domain $D$, every biharmonic mapping
has the representation
$$F(z)=G_1(z)+|z|^2G_2(z),$$ where both $G_1$ and $G_2$ are harmonic in $D$.

More generally, a complex-valued mapping $F$ of a domain $D$ is
called {\it polyharmonic} (or {\it $p$-harmonic}) if $F$ satisfies the
equation $\Delta^{p}F =\Delta(\Delta^{p-1}F)= 0$ for some $p\in \mathbb{N}^{+}$.
In a simply connected domain, a mapping $F$ is polyharmonic if and only if $F$ has
the following representation:
$$F(z) =\sum_{k=1}^{p}|z|^{2(k-1)}G_{k}(z),$$
where each $G_{k}$ is harmonic, i.e., $\Delta G_{k}(z) = 0$ for $k\in \{1,\cdots,p\}$ (see \cite{sh1, sh2}).
Obviously, when $p$=1 (resp. $p$=2), $F$ is harmonic (resp. biharmonic).
The properties of biharmonic mappings have been investigated by many authors (see \cite{z1, z2, z3, CW, jh, s, w}). We refer to
\cite{sh1, sh2, CRW1, CRW2, CRW3} for the discussions on polyharmonic mappings and \cite{jt, du}
for the basic
properties of harmonic mappings.

We use $\mathcal{S}_{H}$ to denote the class consisting of univalent harmonic mappings in $\mathbb{D}$. Such mappings can be
written in the form
\be\label{eq1.2}f(z)=h(z)+\overline{g(z)}=z+\sum_{j=2}^{\infty}a_j
z^j+\sum_{j=1}^{\infty}\overline{b_jz^j},\ee with $|b_{1}|<1$. Let $\mathcal{S}^{\ast}_{H}$ and $\mathcal{C}_{H}$
be the subclasses of $\mathcal{S}_{H}$, where the images of
$f(\mathbb{D})$ are starlike and convex, respectively. If $b_{1}=0$, then $\mathcal{S}_{H}$, $\mathcal{S}^{\ast}_{H}$ and $\mathcal{C}_{H}$ reduce to the classes
$\mathcal{S}_{H}^{0}$, $\mathcal{S}^{0,\ast}_{H}$ and $\mathcal{C}^{0}_{H}$, respectively. See also \cite{du}.

A classical theorem of Fekete and Szeg\H{o} \cite{fe} states that for $f\in \mathcal{S}$ of
the form \eqref{eq1.1}, the functional $|a_{3}-\lambda a_{2}^{2}|$ satisfies the following inequality:
$$|a_{3}-\lambda a_{2}^{2}|\leq \begin{cases}
\displaystyle 3-4\lambda,\;\;\;\;\;\;\;\;\;\;\;\;\;\;\;\;\;\;\;\;\;\lambda\leq0,\\
\displaystyle 1+2e^{-\frac{2\lambda}{1-\lambda}}, \;\;\;\;\;\;\;\; 0\leq \lambda\leq1, \\
\displaystyle 4\lambda-3,\;\;\;\;\;\;\;\;\;\;\;\;\;\;\;\;\;\;\;\;\;\lambda\geq1.
\end{cases}
$$
This inequality is sharp in the sense that for each real $\lambda$ there exists a
function in $\mathcal{S}$ such that equality holds (see \cite{ab, ko1}). Thus the determination
of sharp upper bounds for the nonlinear functional $|a_{3}-\lambda a_{2}^{2}|$ for any
compact family $\mathcal{F}$ of functions in $\mathcal{A}$ is often called the Fekete-Szeg\H{o}
problem for $\mathcal{F}$. Many researchers have studied the Fekete-Szeg\H{o}
problem for analytic close-to-convex functions (see \cite{km, ko1, ko2}).
A natural question is whether we can get similar generalizations
to harmonic close-to-convex mappings.

In \cite{jt}, Clunie and Sheil-Small obtained the following result:

\begin{prop}\label{prop1} {\rm \bf (\cite[Lemma 5.11]{jt})}
If $f=h+\overline{g}\in \mathcal{C}_{H}$, then there exist angles $\alpha$ and $\beta$ such that
$$\mathrm{Re}\left\{\left(e^{i\alpha}h'(z)+e^{-i\alpha}g'(z)\right)(e^{i\beta}-e^{-i\beta}z^{2})\right\}>0$$
for all $z\in \mathbb{D}$.
\end{prop}

Our purpose of this paper is twofold. In Section \ref{csw-sec3}, we obtain coefficient estimates
of a class of close-to-convex harmonic mappings
and, as an application, show upper bounds for the Fekete-Szeg\H{o} functionals $|a_{3}-\lambda a_{2}^{2}|$ and $|b_{3}-\lambda b_{2}^{2}|$.
The main results in this section are
Theorems \ref{thm1} and \ref{thm2}. In Section \ref{csw-sec4},
first we obtain two sufficient conditions for mappings $F\in \mathcal{H}_{p}$ to be starlike with respect to the origin and convex, respectively,
given in Theorems are \ref{thm3} and \ref{thm4}.
Then, we establish some coefficient estimates for two classes of polyharmonic mappings $\mathcal{HS}_{p}$ and $\mathcal{HC}_{p}$.
Our main results are
Theorems \ref{thm5} and \ref{thm6}. Finally, we obtain a  generalization of Proposition \ref{prop1} to the class $\mathcal{HC}_{p}$, which is
Theorem \ref{thm7}.

\section{Preliminaries}\label{csw-sec2}

In \cite{av}, Avci and Z{\l}otkiewicz introduced the class $\mathcal{HS}$ of univalent harmonic mappings $f$ with the series expansion
\eqref{eq1.2} such that
$$\sum_{j=2}^{\infty}j(|a_{j}|+|b_{j}|)\leq 1-|b_{1}|,\quad (0\leq|b_{1}|<1),$$
and the subclass $\mathcal{HC}$ of $\mathcal{HS}$, where
$$\sum_{j=2}^{\infty}j^{2}(|a_{j}|+|b_{j}|)\leq 1-|b_{1}|,\quad (0\leq|b_{1}|<1).$$ These two classes constitute a harmonic
counterparts of classes introduced by Goodman \cite{go}. They are useful in studying questions of so-called $\delta$-neighborhoods originally considered by Ruscheweyh \cite{ru} (see also \cite{qiwa})
and in constructing explicit $k$-quasiconformal extensions of mappings(see Fait et al. \cite{fa}).

We denote by $\mathcal{H}_{p}$ the set of polyharmonic mappings $F$ in $\mathbb{D}$ with the form:
\be\label{eq2.1}
F(z)=\sum_{k=1}^{p}|z|^{2(k-1)}\big(h_{k}(z)+\overline{g_{k}(z)}\big)
=\sum_{k=1}^{p}|z|^{2(k-1)}\sum_{j=1}^{\infty}(a_{k,j}z^{j}+\overline{b_{k,j}}\overline{z^{j}}), \ee
where $a_{1,1}=1$, $|b_{1,1}|<1$. We say that a univalent polyharmonic mapping $F$ with $F(0) = 0$ is
starlike with respect to the origin if the curve $F(re^{i\theta})$
is starlike with respect to the
origin for each $r\in(0,1)$.
The following result gives a convenient characterization of this property.

\begin{prop}{\rm \bf (\cite{po})}\label{pro2}
If $F$ is univalent, $F(0) = 0$ and
$\frac{\partial}{\partial\theta}\big(\arg F(re^{i\theta})\big) >0$ for
$z = re^{i\theta}\not= 0$, then $F$ is starlike with respect to the origin.
\end{prop}

A univalent
polyharmonic mapping $F$ with $F(0) = 0$ and $\frac{\partial}{\partial \theta}F(re^{i\theta})$ $
\neq0$ whenever $r\in(0,1)$, is said to be convex
if the curve $F(re^{i\theta})$ is convex for each $r\in(0,1)$.

\begin{prop}{\rm \bf (\cite{po})} \label{pro3}
 If $F$ is univalent,
$F(0) = 0$, $\frac{\partial}{\partial \theta}F(re^{i\theta})
\neq0$ whenever $r\in(0,1)$, and $\frac{\partial}{\partial \theta}
 \left[ \arg \left(\frac{\partial}{\partial \theta}F(re^{i\theta})\right) \right] >0$ for
$z = re^{i\theta}\not= 0$, then $F$ is convex.
\end{prop}

In \cite{qiwa}, J. Qiao and X. Wang introduced the class $\mathcal{HS}_{p} $ of
polyharmonic mappings $F$ of the form \eqref{eq2.1} satisfying the condition

\be\label{eq2.2}\begin{cases}
\displaystyle \sum_{k=1}^{p}\sum_{j=2}^{\infty}
\big(2(k-1)+j\big)\big(|a_{k,j}|+|b_{k,j}|\big)\leq1-|b_{1,1}|-\sum_{k=2}^{p}(2k-1)\big(|a_{k,1}|+|b_{k,1}|\big),\\
\displaystyle 0\leq|b_{1,1}|+\sum_{k=2}^{p}(2k-1)\big(|a_{k,1}|+|b_{k,1}|\big)<1,
\end{cases} \ee
and the subclass $\mathcal{HC}_{p}$ of $\mathcal{HS}_{p} $, where
\be\label{eq2.3}
\begin{cases}
\displaystyle \sum_{k=1}^{p}\sum_{j=2}^{\infty}
\big(2(k-1)+j^{2}\big)\big(|a_{k,j}|+|b_{k,j}|\big)\leq1-|b_{1,1}|-\sum_{k=2}^{p}(2k-1)\big(|a_{k,1}|+|b_{k,1}|\big),\\
\displaystyle 0\leq|b_{1,1}|+\sum_{k=2}^{p}(2k-1)\big(|a_{k,1}|+|b_{k,1}|\big)<1.
\end{cases}
\ee

Obviously, for any $F\in \mathcal{HS}_{p}$, we have $|F(z)|<2|z|$ for $z\in \mathbb{D}$.
For $p = 1$, the classes $\mathcal{HS}_{p}$ and $\mathcal{HC}_{p}$ reduce to $\mathcal{HS}$ and
$\mathcal{HC}$, respectively. An important property of these classes is given by the following result.

\begin{Thm}{\rm \bf (\cite[Theorem $3.1$]{qiwa})} \label{thmA}
Suppose $F\in \mathcal{HS}_{p}$. Then $F$ is
univalent and sense preserving in $\mathbb{D}$.
\end{Thm}

\section{Coefficient estimates for a class of close-to-convex harmonic mappings}\label{csw-sec3}

In this section, we will consider the coefficient estimates and the Fekete-Szeg\H{o}
problem of mappings of the class $\mathcal{F}$, defined as follows.
In \cite{bp}, Bharanedhar and Ponnusamy obtained the following result:

\begin{Thm}\label{ThmB} {\rm \bf (\cite[Theorem 1]{bp})}
 Let $f=h+\overline{g}$ be a harmonic mapping of $\mathbb{D}$, with $h'(0)\not=0$, which satisfies
 \be\label{eq3.1}
 g'(z)=e^{i\theta}zh'(z)\;and\;\mathrm{Re}\left(1+z\frac{h''(z)}{h'(z)}\right)>-\frac{1}{2}\ee
for all $z\in \mathbb{D}.$ Then $f$ is a univalent close-to-convex mapping in $\mathbb{D}$.
\end{Thm}

We use $\mathcal{F}$ to denote the class of harmonic mapping $f$ in $\mathbb{D}$ of
the form \eqref{eq1.2}, satisfying \eqref{eq3.1}.
Let $\mathcal{H}$ and $\mathcal{G}$ be the subclasses of $\mathcal{F}$,
where $$\mathcal{H}=\{F=h+\overline{g}:F\in \mathcal{F}\;\; \text{and}\;\; g\equiv0\}$$
and
$$\mathcal{G}=\{F=h+\overline{g}:F\in \mathcal{F}\;\; \text{and}\;\; h\equiv0\}.$$

\begin{lem}\label{lem1}
The classes $\mathcal{H}$, $\mathcal{G}$ and $\mathcal{F}$ are compact.
\end{lem}
\bpf
Suppose that $f_{n}=h_{n}+\overline{g_{n}}\in \mathcal{F}$ and that $f_{n}\rightarrow f$
uniformly on compact subsets of $\mathbb{D}$. It follows from
Hurwitz's theorem that $f$ is harmonic, and therefore it has a canonical representation $f=h+\overline{g}$.
It is easy to see that $h_{n}\rightarrow h$ and $g_{n}\rightarrow g$ locally uniformly,
and that $h'(0)=1$ and $g'(z)=e^{i\theta}zh'(z)$.
Because $$\text{Re}\left(1+z\frac{h_{n}''(z)}{h_{n}'(z)}\right)>-\frac{1}{2} $$
in $\mathbb{D}$, it follows that when $n\rightarrow +\infty$, then we have
\be\label{eq03.1} \text{Re}\left(\frac{3}{2}+z\frac{h''(z)}{h'(z)}\right)\geq0\ee
for all $z\in\mathbb{D}$.
It suffices to show that the equation in \eqref{eq03.1} cannot hold for any $z\in \mathbb{D}$.
Obviously, the function $\text{Re}\left(\frac{3}{2}+z\frac{h''(z)}{h'(z)}\right)$ is harmonic in $\mathbb{D}$.
By the maximum principle,
if $\text{Re}\left(\frac{3}{2}+z_{0}\frac{h''(z_{0})}{h'(z_{0})}\right)=0$ for some $z_{0}\in\mathbb{D}$, then $\text{Re}\left(\frac{3}{2}+z\frac{h''(z)}{h'(z)}\right)\equiv0$, and hence
$$\frac{3}{2}+z\frac{h''(z)}{h'(z)}\equiv i C$$
for some real constant $C$. That is a contradiction. Hence,
$$\text{Re}\left(\frac{1}{2}+z\frac{h''(z)}{h'(z)}\right)>-\frac{1}{2}$$
for all $z\in \mathbb{D}$,
and the proof is complete.
\epf

\bthm\label{thm1}
Let $f$ be of the form \eqref{eq1.2} satisfying \eqref{eq3.1}. Then
\be\label{eq3.3} |a_{j}|\leq \frac{j+1}{2}\;\;\text{and}\;\;|b_{j}|\leq \frac{j-1}{2}\ee
for all $j=1,2,\ldots.$
\ethm
\bpf
Since $\text{Re}\left(\frac{3}{2}+z\frac{h''(z)}{h'(z)}\right)>0$, then there exists an analytic
function $p_{1}(z)=c_{0}+c_{1}z+c_{2}z^{2}+\ldots$, such that
\be\label{eq3.4} z\frac{h''(z)}{h'(z)}=p_{1}(z)-\frac{3}{2}=c_{1}z+c_{2}z^{2}+\ldots,\ee
and $\text{Re}\{p_{1}(z)\}>0$. Then \eqref{eq3.4} implies that
$$j(j+1)a_{j+1}z^{j}=\big(ja_{j}c_{1}+(j-1)a_{j-1}c_{2}+\ldots+a_{1}c_{j}\big)z^{j}$$
for $j=1,2,\ldots,$ and hence,
$$a_{j+1}=\frac{1}{j(j+1)}\sum_{\gamma=1}^{j}\gamma a_{\gamma}c_{j+1-\gamma}.$$
Because $p_{1}(z)=c_{0}+c_{1}z+c_{2}z^{2}+\ldots$, and $\text{Re}$ $p_{1}(z)>0$, then by \cite[Lemma 1, p. 50]{du}, we
have $|c_{j}|\leq2\text{Re}\{c_{0}\}=3$ for all $j=1,2,\ldots$. If we write $a_{1}=1$, it follows that
$$|a_{j}|\leq\frac{j+1}{2}\;\;\;\text{for all}\;j=1,2,\ldots.$$
By \eqref{eq3.1}, $g'(z)=e^{i\theta}zh'(z)$, and therefore $\sum_{j=1}^{\infty}jb_{j}z^{j-1}=e^{i\theta}\sum_{j=1}^{\infty}ja_{j}z^{j}$.
Thus, we have $$b_{1}=0,\;\;\;jb_{j}=e^{i\theta}(j-1)a_{j-1}\;\text{for all}\;j=2,3,\ldots.$$
Then, we obtain $$|b_{j}|=\frac{j-1}{j}|a_{j-1}|\leq\frac{j-1}{2}\;\;\;\text{for all}\;j=1,2,\ldots.$$
\epf

Now, we are ready to establish upper bounds for the Fekete-Szeg\H{o} functionals
$|a_{3}-\lambda a_{2}^{2}|$ and $|b_{3}-\lambda b_{2}^{2}|$.

\bthm\label{thm2}
Let $f$ be of the form \eqref{eq1.2} and satisfy \eqref{eq3.1}. Then
\be\label{eq2}|a_{3}-\lambda a_{2}^{2}|\leq \max\left\{\frac{1}{2},\frac{|8-9\lambda|}{4}\right\}
\;\;\;\;\text{and}\;\;\;\;|b_{3}-\lambda b_{2}^{2}|\leq 1+\frac{|\lambda|}{4}\ee
for all $\lambda\in \mathbb{R}$.
\ethm
\bpf
Let $$p_{2}(z)=\frac{2}{3}\left(\frac{3}{2}+z\frac{h''(z)}{h'(z)}\right),$$
where $h(z)=z+\sum_{j=1}^{\infty}a_{j}z^{j}$. Then by simple calculations, we obtain

\be\label{eq3.5}p_{2}(z)=1+\frac{2}{3}\big(2a_{2}z+(6a_{3}-4a_{2}^{2})z^{2}+\ldots\big).\ee
Write $p_{2}(z)=1+u_{1}z+u_{2}z^{2}+\ldots$. Because $\text{Re}$ $p_{2}(z)>0$, then by \cite[formula (10), p. 166]{po}, we have
$$\left|u_{2}-\frac{ u_{1} ^{2}}{2}\right|\leq 2-\frac{|u_{1}|^{2}}{2}.$$
It follows from \eqref{eq3.5} that $$a_{2}=\frac{3}{4}u_{1}\;\;\text{and}\;\;
a_{3}=\frac{1}{4}u_{2}+\frac{3}{8}u_{1}^{2}.$$
Hence,
\begin{align}\begin{split}\label{eq3.6}
|a_{3}-\lambda a_{2}^{2}|&=\left|\frac{1}{4}u_{2}+\frac{3}{8}u_{1}^{2}-\frac{9}{16}\lambda u_{1}^{2} \right|\\
&=\frac{1}{4}\left|u_{2}-\frac{1}{2}u_{1}^{2} +\frac{8-9\lambda}{4} u_{1}^{2} \right|\\
&\leq \frac{1}{4}\left(2-\frac{1}{2}|u_{1}|^{2} +\frac{|8-9\lambda|}{4}| u_{1}|^{2} \right).\\
\end{split}\end{align}
If $\frac{|8-9\lambda|}{4}<\frac{1}{2}$,
then \eqref{eq3.6} implies $$|a_{3}-\lambda a_{2}^{2}|\leq \frac{1}{2}.$$
Equality is attained if we choose $a_{2}=0$ and $a_{3}=\pm \frac{1}{2}$.

If $\frac{|8-9\lambda|}{4}\geq \frac{1}{2}$, then it follows from \cite[Lemma 1, p. 50]{du} that
$|u_{1}|\leq 2\text{Re}\{p_{2}(0)\}=2$ and \eqref{eq3.6} that
$$|a_{3}-\lambda a_{2}^{2}|
\leq \frac{1}{2}+\frac{1}{4}\left(\frac{|8-9\lambda|}
{4}-\frac{1}{2}\right)|u_{1}|^{2}\leq\frac{|8-9\lambda|}{4}.$$
Choosing $a_{2}=\pm \frac{3}{2}$ and $a_{3}=2$ in \eqref{eq3.6} shows that the result is sharp.

Since $g'(z)=e^{i\theta}zh'(z)$,
we have $$\sum_{j=1}^{\infty}jb_{j}z^{j-1}= e^{i\theta}\sum_{j=1}^{\infty}ja_{j}z^{j}.$$
Obviously, $b_{2}=\frac{ e^{i\theta}}{2}a_{1}=\frac{ e^{i\theta}}{2}$ and $b_{3}=\frac{2 e^{i\theta}}{3}a_{2}$.
Hence, \eqref{eq3.3} implies
$$|b_{3}-\lambda b_{2}^{2}|=\left|\frac{2 e^{i\theta}}{3}a_{2}-\frac{\lambda e^{2i\theta}}{4} \right|\leq \frac{2}{3}|a_{2}|+\frac{|\lambda|}{4}\leq1+\frac{|\lambda|}{4}.$$
If $\lambda\geq0$, then equality is attained when $b_{3}=-e^{2i\theta}$, i.e. $a_{2}=-\frac{3}{2}e^{i\theta}$. If $\lambda<0$, then equality is attained when $b_{3}=e^{2i\theta}$, i.e. $a_{2}=\frac{3}{2}e^{i\theta}$.
\epf

\begin{rem}\label{rem1}
Both equalities in \eqref{eq2} are attained when $a_{2}=\frac{3}{2}$ and $b_{3}=e^{i\theta}$ or $a_{2}=-\frac{3}{2}$ and $b_{3}=-e^{i\theta}$,
but only in the case $|8-9\lambda|\geq 2$ and $\theta=2k\pi$, where $k\in \mathbb{Z}$.
\end{rem}

\section{Coefficient estimates for two classes of polyharmonic mappings}\label{csw-sec4}

Let $L$ denote the following differential operator defined on the class of complex-valued $C^{1}$ functions:
$$L=z\frac{\partial}{\partial z}-\overline{z}\frac{\partial}{\partial \overline{z}}.$$

An important property of the operator $L[F]$ is that it
 behaves with respect to polyharmonic mappings much like the operator $zf'(z)$ defined for analytic functions
 (see {\rm \bf \cite[Corollary 1(3)]{z3}}).

\begin{thm}\label{thm3} Each mapping $F \in \mathcal{HS} _{p}$ is starlike with respect to the origin.
\end{thm}
\bpf
Let $F \in \mathcal{HS} _{p}$ be of the form \eqref{eq2.1}. It follows from Theorem A that
\be\label{eq4.1}J_{F}(0)=1-|b_{1,1}|^{2}>0.\ee
By computation, we have
\begin{align}\begin{split}\label{eq4.2}
\frac{\partial}{\partial\theta}\big(\arg F (re^{i\theta})\big)
&=\frac{\partial}{\partial\theta}\left[\text{Im} \big(\log F (re^{i\theta})\big)\right]\\
&=\text{Im}\left[\frac{\partial}{\partial\theta} \big(\log F (re^{i\theta})\big)\right]\\
&=\text{Im}\frac{izF_{z}(z)-i\overline{z}F_{\overline{z}}(z)}{F(z)}\\
&=\text{Re} \frac{L[F (z)]}{ F (z)}\\
&=\text{Re} \frac{\sum_{k=1}^{p}|z|^{2(k-1)}\sum_{j=1}^{\infty}\big(ja_{k,j}z^{j}-j\overline{b_{k,j}z^{j}}\big)}
{\sum_{k=1}^{p}|z|^{2(k-1)}\sum_{j=1}^{\infty}
\big( a_{k,j}z^{j}+\overline{b_{k,j}z^{j}}\big)}.\\
\end{split}\end{align}
By \eqref{eq4.1} and \eqref{eq4.2}, we obtain
\begin{align}\begin{split}\label{eq4.3}
\lim_{z\rightarrow0}\frac{\partial}{\partial\theta}\big(\arg F (re^{i\theta})\big)
&=\lim_{z\rightarrow0}\text{Re} \frac{1-\overline{b_{1,1}}\overline{z}/z}{ 1+\overline{b_{1,1}}\overline{z}/z}\\
&=\lim_{z\rightarrow0}\text{Re} \frac{\left(1-\overline{b_{1,1}}\overline{z}/z\right)\left(1+b_{1,1}z/\overline{z}\right)}{ \left|1+\overline{b_{1,1}}\overline{z}/z\right|^2}\\
&=\lim_{z\rightarrow0}\frac{1-|b_{1,1}|^2}{| 1+\overline{b_{1,1}}\overline{z}/z|^2}\\
&\geq\frac{1-|b_{1,1}|^{2}}{(1+|b_{1,1}|)^{2}}>0.\\
\end{split}\end{align}
It follows from Theorem A that each $F \in \mathcal{HS} _{p}$ is univalent in $\mathbb{D}$.
Then, we have that $F(z)\neq 0$ for $z\in \mathbb{D}\setminus \{0\}$, and the function $\text{Re} \frac{L[F (z)]}{ F (z)}$ is continuous in $\mathbb{D}\setminus\{0\}$.
Therefore, by \eqref{eq4.2}, \eqref{eq4.3} and the continuity of $F$ in $\mathbb{D}\setminus\{0\}$, we see the condition $\frac{\partial}{\partial\theta}\big(\arg F(re^{i\theta})\big)>0$ for all $z\in \mathbb{D}\setminus\{0\}$ is equivalent to
\be\label{eq4.4} \frac{L[F(z)]}{ F(z)}\neq \frac{\zeta-1}{\zeta+1}\ee
for all $z\in \mathbb{D}\setminus\{0\}$ and all $\zeta\in \mathbb{C}$ with $|\zeta|=1$ and $\zeta\neq-1$. Hence, \eqref{eq4.4} holds if and only if
$$\Phi(z):=(\zeta+1)L[F(z)]-(\zeta-1)F(z)\neq 0$$
for all $z\in \mathbb{D}\setminus\{0\}$ and all $|\zeta|=1$. Calculations show
that
\begin{align*}
|\Phi(z)|
=&\left|(\zeta+1)\sum_{k=1}^{p}|z|^{2(k-1)}\sum_{j=1}^{\infty} j \big(a_{k,j}z^{j}-\overline{b_{k,j}z^{j}}\big)\right.\\
&\left.- (\zeta-1)\sum_{k=1}^{p}|z|^{2(k-1)}\sum_{j=1}^{\infty}  \big(a_{k,j}z^{j}+\overline{b_{k,j}z^{j}}\big)\right|\\
=&\left|
\sum_{k=1}^{p}|z|^{2(k-1)}\sum_{j=1}^{\infty}\left(\big(j+1+\zeta(j-1)\big)a_{k,j}z^{j}
-\big(j-1+\zeta(j+1)\big)\overline{b_{k,j}z^{j}} \right)
\right|.\\
\end{align*}
If $F$ is the identity, obviously, we have $|\Phi(z)|=2|z|$.
If $F(z)=z+\overline{b_{1,1}z}$, then $$|\Phi(z)|=2\left|z-\zeta\overline{b_{1,1}z}\right|\geq2|z|(1-|b_{1,1}|).$$
If $F$ is not an affine mapping, then
$$|\Phi(z)|>|z|\left(2-2|b_{1,1}|-2\sum_{k=1}^{p}\sum_{j=2}^{\infty}j(|a_{k,j}|+|b_{k,j}|)
-2\sum_{k=2}^{p}(|a_{k,1}|+|b_{k,1}|) \right). $$
Hence, each mapping $F\in \mathcal{HS}_{p}$ is starlike with respect to the point $z=0$.
\epf

\begin{thm}\label{thm4} Each mapping $F \in \mathcal{HC} _{p}$ is convex.
\end{thm}
\bpf
Let $F \in \mathcal{HC }_{p}$ be of the form \eqref{eq2.1}. By \eqref{eq4.2}, we have
\begin{align}\begin{split}\label{eq4.5}
\frac{\partial}{\partial \theta}
\left[ \arg \left(\frac{\partial}{\partial \theta}F(re^{i\theta})\right) \right]
&=\text{Re} \frac{L\left[\frac{\partial}{\partial\theta}F(re^{i\theta})\right ]}{\frac{\partial}{\partial\theta}F(re^{i\theta})}\\
&=\text{Re} \frac{L[L[F(z)]]}{L[F(z)]}\\
&=\text{Re} \frac{\sum_{k=1}^{p}|z|^{2(k-1)}\sum_{j=1}^{\infty}\big(j^{2}a_{k,j}z^{j}+j^{2}\overline{b_{k,j}z^{j}}\big)}
{\sum_{k=1}^{p}|z|^{2(k-1)}\sum_{j=1}^{\infty}
\big(ja_{k,j}z^{j}-j\overline{b_{k,j}z^{j}}\big)}.\\
\end{split}\end{align}
Then by \eqref{eq4.1} and \eqref{eq4.5}, we have
\begin{align}\begin{split}\label{eq4.6}
\lim_{z\rightarrow0}\frac{\partial}{\partial \theta}
\left[ \arg \left(\frac{\partial}{\partial \theta}F(re^{i\theta})\right) \right]
&=\lim_{z\rightarrow0}\text{Re} \frac{1+\overline{b_{1,1}}\overline{z}/z}{ 1-\overline{b_{1,1}}\overline{z}/z}\\
&=\lim_{z\rightarrow0}\text{Re} \frac{\left(1+\overline{b_{1,1}}\overline{z}/z\right)\left(1-b_{1,1}z/\overline{z}\right)}{ \left|1-\overline{b_{1,1}}\overline{z}/z\right|^2}\\
&=\lim_{z\rightarrow0}\frac{1-|b_{1,1}|^2}{| 1-\overline{b_{1,1}}\overline{z}/z|^2}\\
&\geq\frac{1-|b_{1,1}|^{2}}{(1+|b_{1,1}|)^{2}}>0.\\
\end{split}\end{align}
If $F$ is the identity, obviously, we have $|L[F(z)]|=|z|$.
If $F(z)=z+\overline{b_{1,1}z}$, then $$|L[F(z)]|=|z-\overline{b_{1,1}z}|\geq|z|(1-|b_{1,1}|).$$
If $F$ is not an affine mapping, then
\begin{align*}
|L[F(z)]|&=\left|\sum_{k=1}^{p}|z|^{2(k-1)}\sum_{j=1}^{\infty}
\big( ja_{k,j}z^{j}-j\overline{b_{k,j}z^{j}}\big)\right|\\
&>|z|\left(1-|b_{1,1}|-\sum_{k=1}^{p}\sum_{j=2}^{\infty}
j\big( |a_{k,j}|+|b_{k,j}|\big)-\sum_{k=2}^{p}
j\big( |a_{k,1}|+|b_{k,1}|\big)\right).
\end{align*}
Therefore, $F\in \mathcal{HC}_{p}$ implies $L[F(z)]\neq0$ for $z\in \mathbb{D}\setminus \{0\}$,
and hence the function $\text{Re} \frac{L[L[F (z)]]}{ L[F (z)]}$ is continuous in $\mathbb{D}\setminus\{0\}$.
Therefore, by \eqref{eq4.5}, \eqref{eq4.6} and the continuity of $L[F(z)]$ in $\mathbb{D}\setminus\{0\}$,
we see the condition $\frac{\partial}{\partial \theta}\left[ \arg \left(\frac{\partial}{\partial \theta}F(re^{i\theta})\right) \right]>0$
for all $z\in \mathbb{D}\setminus\{0\}$ is equivalent to
\be\label{eq4.7}\text{Re} \frac{L[L[F(z)]]}{L[F(z)]}\neq \frac{\zeta-1}{\zeta+1}\ee
for all $z\in \mathbb{D}\setminus\{0\}$ and all $\zeta\in \mathbb{C}$ with $|\zeta|=1$ and $\zeta\neq-1$.
Hence, \eqref{eq4.7} holds if and only if
$$\Psi(z):=(\zeta+1)L[L[F(z)]]-(\zeta-1)L[F(z)]\neq 0$$
for all $z\in \mathbb{D}\setminus\{0\}$ and all $|\zeta|=1$. Calculations show
that
\begin{align*}
|\Psi(z)|
=&\left|(\zeta+1)\sum_{k=1}^{p}|z|^{2(k-1)}\sum_{j=1}^{\infty} j^{2}\big(a_{k,j}z^{j}+\overline{b_{k,j}z^{j}}\big)\right.\\
&\left.- (\zeta-1)\sum_{k=1}^{p}|z|^{2(k-1)}\sum_{j=1}^{\infty} j\big(a_{k,j}z^{j}-\overline{b_{k,j}z^{j}}\big)\right|\\
=&\sum_{k=1}^{p}|z|^{2(k-1)}\sum_{j=1}^{\infty} \Big(\big(j^{2}+j+\zeta(j^{2}-j)\big)a_{k,j}z^{j}
+\big(j^{2}-j+\zeta(j^{2}+j)\big)\overline{b_{k,j}z^{j}}\Big).\\
\end{align*}
If $F$ is the identity, obviously, we have $|\Psi(z)|=2|z|$.
If $F(z)=z+\overline{b_{1,1}z}$, then $$|\Psi(z)|=2\left|z+\zeta\overline{b_{1,1}z}\right|\geq2|z|(1-|b_{1,1}|).$$
If $F$ is not an affine mapping, then
$$|\Psi(z)|>|z|\left(2-2|b_{1,1}|-2\sum_{k=1}^{p}\sum_{j=2}^{\infty}j^{2}
\big(|a_{k,j}|+|b_{k,j}|\big)-2\sum_{k=2}^{p}\big(|a_{k,1}|+|b_{k,1}|\big) \right). $$
It follows that each mapping $F\in \mathcal{HC}_{p}$ is convex.
\epf

\begin{figure}
\centering
\includegraphics[width=2.5in]{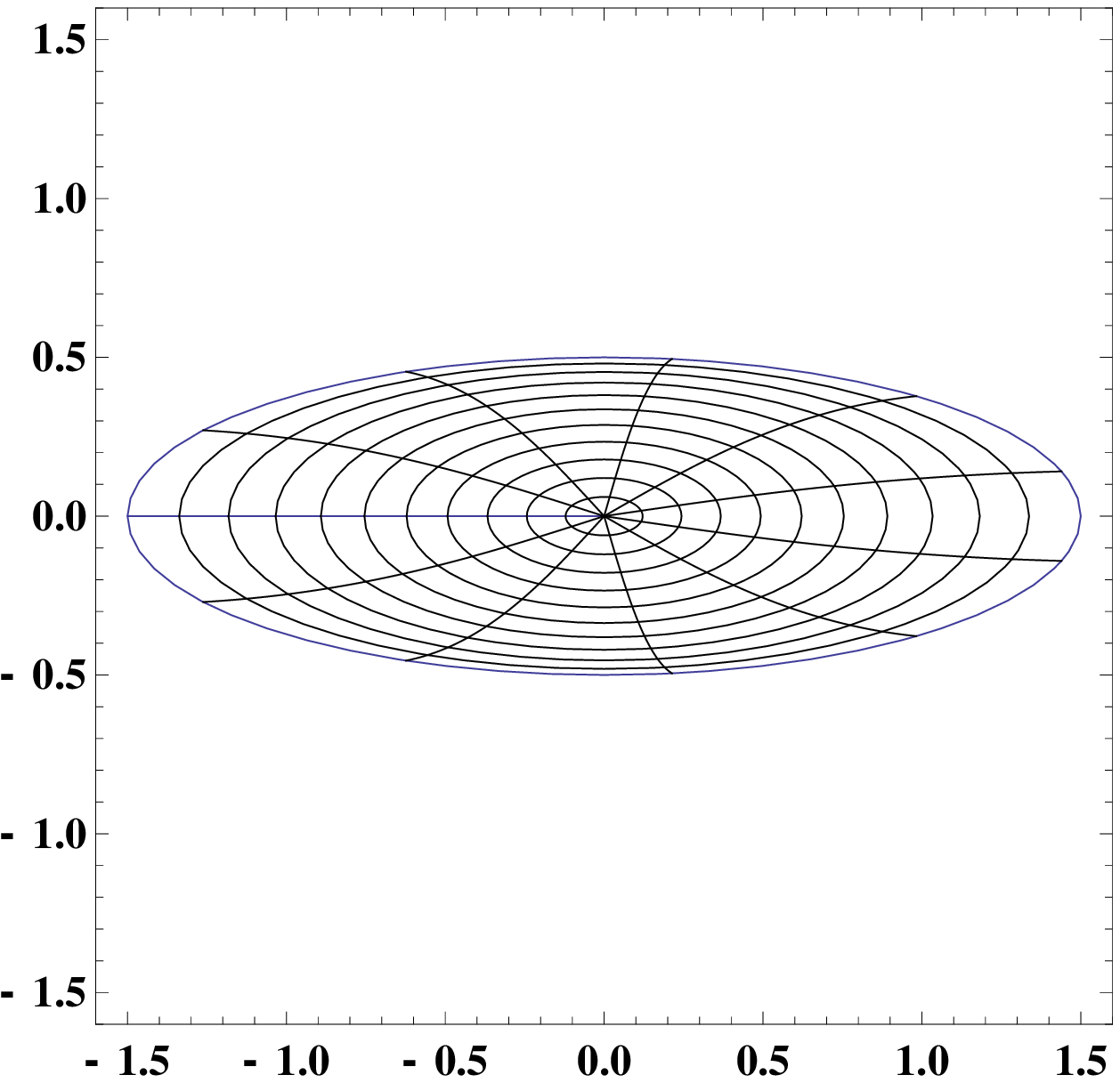}
\quad
\includegraphics[width=2.5in]{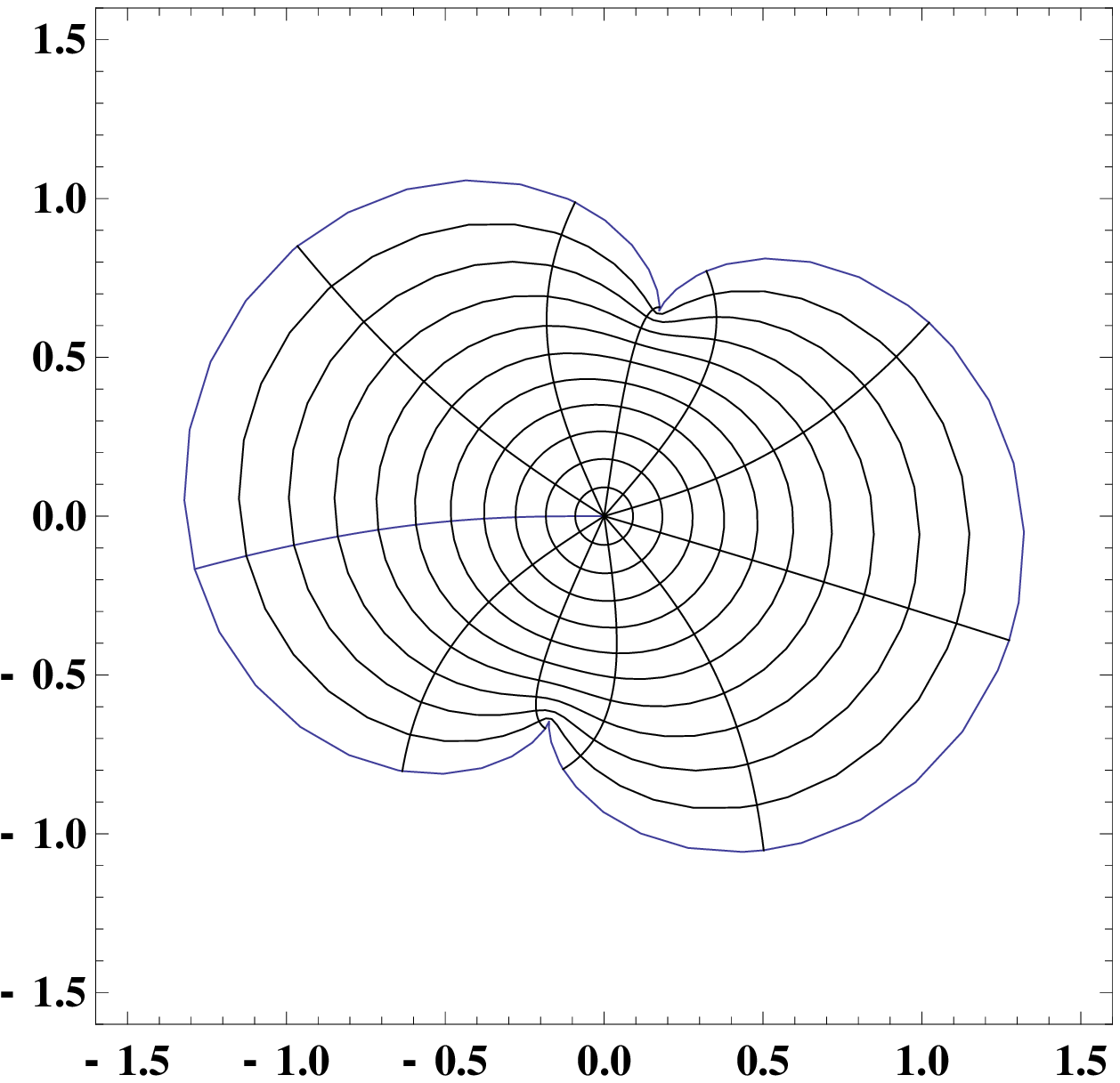}
\caption{
The images of $\mathbb{D}$ under the mappings $F_{1}(z)=z+\frac{1}{3}\overline{z}+\frac{1}{6}|z|^2\overline{z}$ (left) and $F_2(z)=z+\frac{z^{3}}{3}e^{\frac{\pi}{6} i}.$
}\label{figure1}
\end{figure}

\begin{examp} Let $F_{1}(z)=z+\frac{1}{3}\overline{z}+\frac{1}{6}|z|^2\overline{z}$. Then $F_{1}$
is convex. See also Figure \ref{figure1}.
\end{examp}

It is well known that the coefficients of every starlike mapping $f\in \mathcal{S}_{H}^{\ast,0}$ of the form \eqref{eq1.2}
satisfy the sharp inequalities
$$|a_{j}|\leq\frac{(2j+1)(j+1)}{6},\;\;\;|b_{j}|\leq\frac{(2j-1)(j-1)}{6},\;\;\; \big||a_{j}|-|b_{j}|\big|\leq j$$
for $j=2,3,\ldots$ (see \cite{shel}).
The coefficients of each mapping $f\in \mathcal{C}_{H}^{0}$ satisfy the sharp inequalities
$$|a_{j}|\leq\frac{j+1}{2},\;\;\;|b_{j}|\leq\frac{j-1}{2},\;\;\;\text{and}\;\;\;\big||a_{j}|-|b_{j}|\big|\leq 1$$
for $j=2,3,\ldots$ (see \cite{jt}).

Next, we obtain similar results for mappings in $\mathcal{HS}_{p}$ and $\mathcal{HC}_{p}$.

\bthm\label{thm5}
The coefficients of every mapping $F\in \mathcal{HS}_{p}$ satisfy the sharp inequalities
\be\label{eq4.8}\sum_{k=1}^{p}(|a_{k,j}|+|b_{k,j}|)\leq \frac{1}{j}  \ee
for all $j=2,3,\ldots$.
\ethm
\bpf
Let $F\in \mathcal{HS}_{p}$ be of the form \eqref{eq2.1}. By \eqref{eq2.2}, we have
\begin{align*}
\sum_{k=1}^{p}j (|a_{k,j}|+|b_{k,j}|)
&\leq  \sum_{k=1}^{p}\sum_{j=2}^{\infty}\big (2(k-1)+j \big)\big(|a_{k,j}|+|b_{k,j}|\big)
\leq1.
\end{align*}
It follows that $$\sum_{k=1}^{p}  (|a_{k,j}|+|b_{k,j}|)\leq \frac{1}{j } $$
for $j=2,3,\ldots$.
\epf

\begin{examp} Let $F_{2}(z)=z+\frac{z^{j}}{j}e^{i\varphi}$ for all $j=2,3,\ldots$ and $\varphi\in \mathbb{R}$.
Then $F_{2}\in \mathcal{HS}$ is univalent, sense
preserving and starlike with respect to the origin. Obviously, the coefficients of $F_{2}$ satisfy \eqref{eq4.8}. See Figure \ref{figure1} for the case where $j=3$ and $\varphi=\pi/6$.
\end{examp}
The above example shows that the coefficient estimate \eqref{eq4.8} is sharp for $p=1$.

\bthm\label{thm6} The coefficients of each mapping $F\in \mathcal{HC}_{p}$ satisfy the sharp inequalities
\be\label{eq4.9} \sum_{k=1}^{p}  (|a_{k,j}|+|b_{k,j}|)\leq \frac{1}{j^{2}}  \ee
for $j=2,3,\ldots$.
\ethm
\bpf
The proof of Theorem \ref{thm6} is similar to the proof of Theorem \ref{thm5}, and we will omit it.
\epf

\begin{figure}
\centering
\includegraphics[width=2.5in]{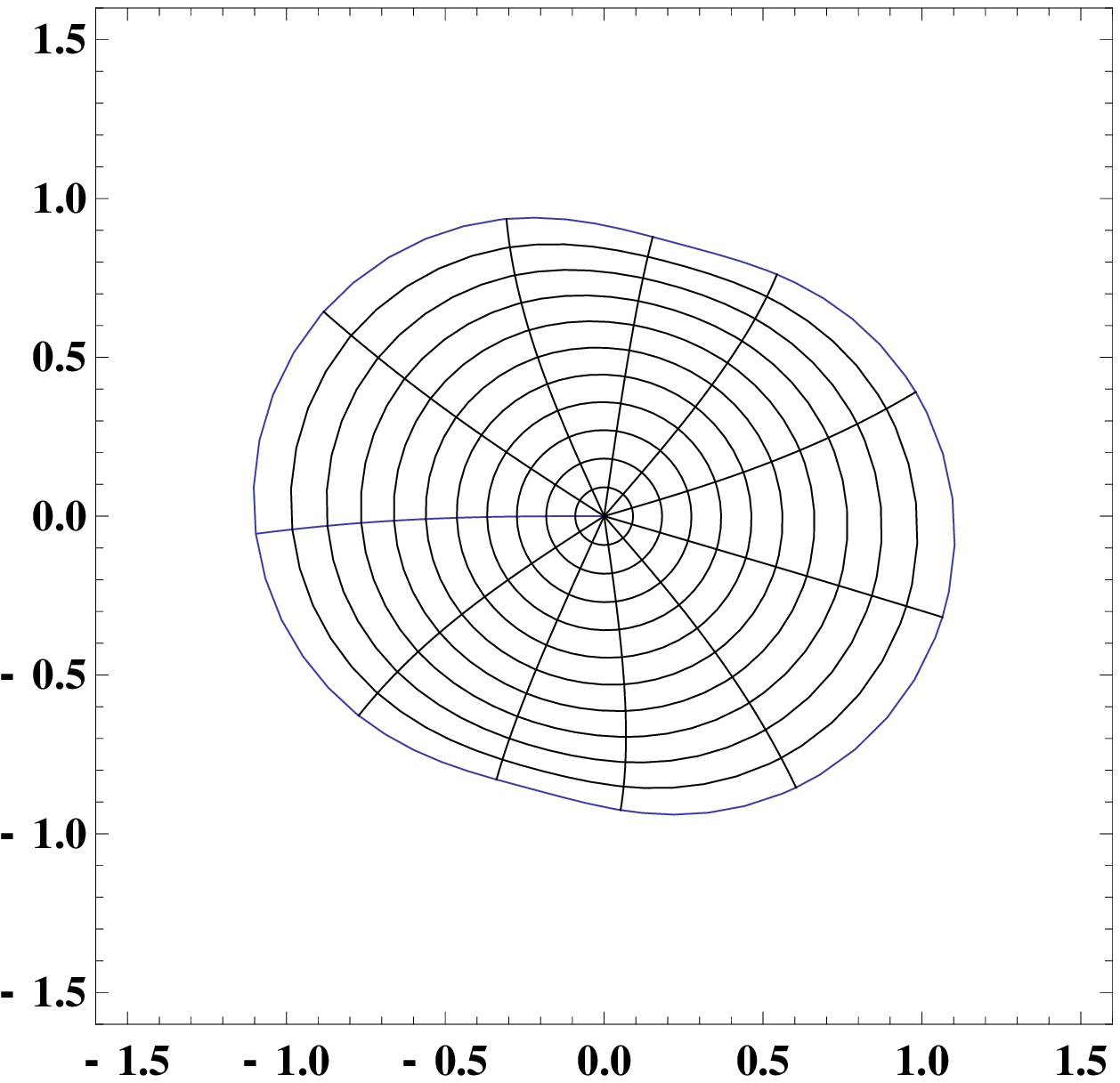}\quad
\includegraphics[width=2.5in]{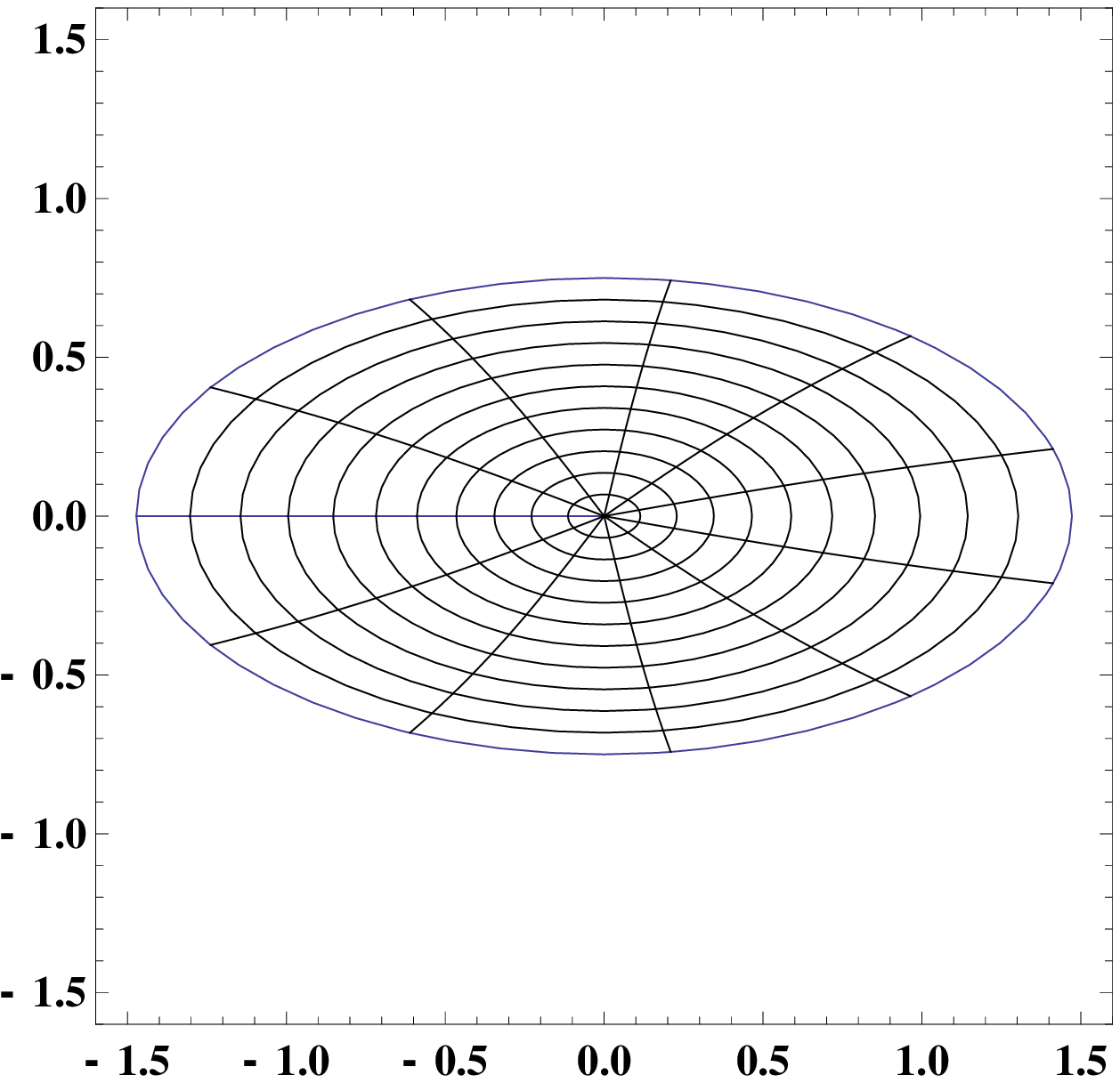}
\caption{The images of $\mathbb{D}$ under the mappings $F_3(z)=z+\frac{z^{3}}{9}e^{ \frac{\pi}{6}i}$ (left) and $F_{4}(z)=z+\frac{1}{9}|z|^{2}z+\frac{1}{4}\overline{z}+\frac{1}{9}|z|^2\overline{z}$.}\label{figure2}
\end{figure}

\begin{examp} Let $F_{3}(z)=z+\frac{z^{j}}{j^{2}}e^{i\varphi}$ for all $j=2,3,\ldots$ and $\varphi\in \mathbb{R}$.
Then $F_{3}\in \mathcal{HC}$ is a univalent, sense
preserving and convex harmonic mapping. Obviously, the coefficients of $F_{3}$ satisfy \eqref{eq4.9}. See Figure \ref{figure2} for the case where $j=3$ and $\varphi=\pi/6$.
\end{examp}
This example shows that the coefficient estimate \eqref{eq4.9} is sharp for $p=1$.

Now, we are ready to generalize Proposition 1 to the polyharmonic mappings of the class $\mathcal{HC}_{p}$.
\bthm\label{thm7}
If $F\in \mathcal{HC}_{p}$ and $a_{k,1}=0$ for $k\in\{2,\ldots,p\}$, then there exist angles $\alpha$ and $\beta$ such that
\be\label{eq4.10}\mathrm{Re}\left\{\left(e^{i\alpha}\sum_{k=1}^{p}|z|^{2(k-1)}h'_{k}(z)+e^{-i\alpha}
\sum_{k=1}^{p}|z|^{2(k-1)}g'_{k}(z)\right)(e^{i\beta}-e^{-i\beta}z^{2})\right\}>0\ee
for all $z\in \mathbb{D}$.
\ethm
\bpf
Let $F\in \mathcal{HC}_{p}$ be of the form \eqref{eq2.1}, fix $r\in(0,1)$, and let
$$F_{r}(z)=\sum_{k=1}^{p}r^{2(k-1)}\big(h_{k}(z)+\overline{g_{k}(z)}\big)=\sum_{j=1}^{\infty}\sum_{k=1}^{p}
\left(a_{k,j}r^{2(k-1)}z^{j}
+\overline{b_{k,j}}r^{2(k-1)}\overline{z^{j}}\right),\; z\in \mathbb{D}.
$$
Then $F_{r}$ is harmonic. By the hypothesis and \eqref{eq2.3}, $F\in \mathcal{HC}_{p}$ implies
$$\sum_{j=2}^{\infty}j^{2}\left|\sum_{k=1}^{p}a_{k,j} r^{2(k-1)}\right|
+\sum_{j=2}^{\infty}j^{2}\left|\sum_{k=1}^{p}b_{k,j}r^{2(k-1)} \right|\leq1-\left|\sum_{k=1}^{p}b_{k,j}r^{2(k-1)} \right|,$$
i.e., $F_{r}\in \mathcal{C}_{H}$ (see \cite{av}). Then Proposition 1 implies that
there exist angles $\alpha$ and $\beta$ such that
$$\text{Re}\left\{\left(e^{i\alpha}\sum_{k=1}^{p}r^{2(k-1)}h'_{k}(z)+e^{-i\alpha}
\sum_{k=1}^{p}r^{2(k-1)}g'_{k}(z)\right)(e^{i\beta}-e^{-i\beta}z^{2})\right\}>0$$
for all $z\in \mathbb{D}$. Let $r=|z|$. The result is proved.
\epf

\begin{examp} Obviously, the mapping $F_{1}(z)=z+\frac{1}{3}\overline{z}+\frac{1}{6}|z|^2\overline{z}\in \mathcal{HC}_{2}$ (see Figure \ref{figure1}). Let $\alpha=\beta=0$. Then $F_{1}$ satisfies the inequality \eqref{eq4.10}.

However, the mapping $F_{4}(z)=z+\frac{1}{9}|z|^{2}z+\frac{1}{4}\overline{z}+\frac{1}{9}|z|^2\overline{z}\in \mathcal{HC}_{2}$ (see Figure \ref{figure2}) also satisfies the inequality \eqref{eq4.10} for $\alpha=\beta=0$ with $a_{2,1}=\frac{1}{9}$.
\end{examp}

\begin{rem}
The proof of Theorem \ref{thm7} requires a somewhat unnatural additional assumption concerning the coefficients $a_{k,1}$. It is not obvious if the result holds without this assumption.
\end{rem}

\end{document}